\input amstex
\documentstyle{amsppt}
\magnification 1200

\NoBlackBoxes
\NoRunningHeads
\topmatter
\title Quenched central limit theorems for a stationary linear process
\endtitle
\author Dalibor Voln\'y and Michael Woodroofe\endauthor
\address{Universit\'e de Rouen, Laboratoire Rapha\"el Salem
UMR 6085 CNRS, Universit\'e de Rouen, Saint \'Etienne-du-Rouvray, France 
\smallskip
University of Michigan, Department of Statistics, Ann Arbor, MI, USA}
\endaddress

\abstract
We establish a sufficient condition under which a central limit theorem for a a stationary linear process is quenched.
We find a stationary linear process for which the Maxwell-Woodroofe's condition is satisfied, $\sigma_n=\|S_n\|_2 = o(\sqrt n)$,
$S_n/\sigma_n$ converge to the standard normal law, and the convergence is not quenched; the weak invariance principle does not hold.
\endabstract

\endtopmatter

\document

\subheading{1. Introduction} 

Let $T$ be an ergodic automorphism of a probability space $(\Omega, \Cal A, \mu)$.  
For $h\in L^2$, $Uh= h\circ T$ is a unitary operator; we will freely switch from the notation $h\circ T^i$ to $U^ih$
and vice versa.

Let $(\Cal F_i)_i$ be
a filtration such that $\Cal F_{i+1} = T^{-1}\Cal F_i$, and $e\in L^2(\Cal F_0)\ominus L^2(\Cal F_{-1})$.
For simplicity we will suppose $\|e\|_2 = 1$.
Let $a_i$ be real numbers with $\sum_{i\in \Bbb N} a_i^2 <\infty$ and let
$$
  f = \sum_{i\leq 0} a_{-i} U^ie.
$$
Then $f\in L^2$ and we say that $(f\circ T^i)_i$ is a causal stationary linear process. 
The stationary linear process is a classical and important case of a (strictly) stationary process and, moreover, any regular stationary 
process is a sum of stationary linear process ``living'' in mutually orthogonal and $U$-invariant subspaces of $L^2$ (cf\. \cite{VWoZ}). 
If $e_k\in  L^2(\Cal F_0)\ominus L^2(\Cal F_{-1})$, $\|e_k\|_2 = 1$ are  mutually orthogonal, 
$a_{k,i}$ are real numbers with $\sum_{k=1}^\infty \sum_{i\in \Bbb N} a_{k,i}^2 <\infty$, and if
$$
  f = \sum_{k=1}^\infty \sum_{i\leq 0} a_{k,-i} U^ie_k \tag{1}
$$
then we say that $(f\circ T^i)_i$ is a causal superlinear process. As shown in \cite{VWoZ}, if $f\in L^2$ is $\Cal F_0$-measurable
and $E(f | \Cal F_{-\infty}) = 0$ (i.e\. the process $(f\circ T^i)_i$ is regular) then a representation (1) exists. 
%Any adapted ($f$ is $\Cal F_0$-measurable) and regular
%($E(f | \Cal F_{-\infty}) = 0$) $L^2$ stationary process $(f\circ T^i)_i$ is thus a sum of mutually orthogonal stationary linear processes.
%Stationary linear processesare thus building blocs of general regular processes.

Let us denote $S_n(f) = \sum_{i=0}^{n-1} f\circ T^i$.
Recall (\cite{PU06}) that if $\sigma_n = \|S_n(f)\|_2 \to \infty$ then the distributions of $S_n(f)/\sigma_n$ weakly
converge to $\Cal N(0,1)$, i.e\. we have a CLT.
We will study when this CLT is quenched. 
%After presenting a sufficient condition we will show an example where the CLT is not quenched
%and also the weak invariance principle does not hold.

Let us suppose that the regular conditional probabilities $m_\omega$ with respect to the $\sigma$-field $\Cal F_0$ exist. If for 
$\mu$ a.e\. $\omega$ the distributions $m_\omega(S_n(f)/\sigma_n)^{-1}$ weakly converge to $\Cal N(0,1)$, we say that the CLT is quenched.
A quenched CLT can be defined using Markov Chains. Any stationary process $(f\circ T^i)_i$ can be expressed using a homogeneous and 
stationary Markov Chain $(\xi_i)_i$ as $(g(\xi_i))_i$; a CLT is quenched if it takes place for a.e\. starting point (this approach 
is probably earlier than our one; it has been used in e.g\. \cite{DLi}). Equivalence of both approaches was explained e.g\. 
in \cite{V}.

In the next section, for a stationary linear processes, we will give a sufficient condition for a quenched CLT. Then, in Section 3,
we will present a stationary linear processes $(f\circ T^i)_i$ for which the Maxwell-Woodroofe's condition and the Hannan's condition
are satisfied but the CLT is not quenched and the weak invariance principle (WIP) does not hold. 
As we will explain in Remark 3, for norming by $\sqrt n$ the limit behaviour of the process is different.
\comment
As shown in \cite{CuMe}, 
Maxwell-Woodroofe's condition implies a quenched CLT and WIP for $S_n(f)/\sqrt n$ (cf\. also \cite{PU}). 
In our example, the Maxwell-Woodroofe's condition is 
satisfied, $S_n(f)/\sqrt n$ converge to zero, $S_n(f)/\sigma_n$ converge to $\Cal N(0,1)$, but the CLT is not quenched and the weak invariance 
principle does not hold. 
\endcomment
\bigskip

\subheading{2. A Sufficient Condition}
Let $(e\circ T^i)$ be a martingale difference sequence as defined in the introduction,
%$e_k = e\circ T^k$ be a sequence of martingale differences adapted to a filtration $(\Cal F_i)_i$, $\Cal F_i = T^{-i}\Cal F_0$, 
$$
  f = \sum_{i=0}^\infty a_ie\circ T^{-i}
$$
where $\sum_{i=0}^\infty  a_i^2 < \infty$, and let $S_n = \sum_{i=0}^{n-1} f\circ T^i$, $n=1,2,\dots$. By definition,
$$\multline
  S_n = \sum_{j=0}^{n-1} f\circ T^j = \sum_{i=0}^\infty \sum_{j=0}^{n-1} a_ie\circ T^{j-i} = 
  \sum_{k=-\infty}^{n-1} \sum_{j=k\vee 0}^{n-1} a_{j-k}e\circ T^k = \\
  = \sum_{k=1}^{n-1} b_{n-k} e\circ T^k + \sum_{k=-\infty}^{0} (b_{n-k}-b_{-k})  e\circ T^k
  \endmultline
$$  
where
$$
  u\vee v = \max\{u,v\},\,\,\, b_0=0,\,\,\, b_j = \sum_{i=0}^{j-1} a_i,\quad j\geq 1.
$$  

We denote 
$$
  \bar{\sigma}_n^2 = E[(S_n - E(S_n| \Cal F_0))^2] = \sum_{k=1}^{n-1} b_{n-k}^2.
$$  

\proclaim{Theorem 1} Let $\bar{\sigma}_n^2 \to \infty$. If
\roster
\item"(i)" $e\circ T^i$ are iid and $(\Cal F_i)_i$ is the natural filtration \newline
or if
\item"(ii)" 
$$
  \sup_{n\geq 1} \max_{k\leq n} \frac{nb_k^2}{\bar{\sigma}_n^2} = c < \infty \tag{2}
$$
\endroster
then for $(1/ \bar{\sigma_n}) [S_n - E(S_n| \Cal F_0)]$ a quenched CLT holds true.
\endproclaim

\underbar{Remark 1}. If the sums $b_k$ converge to a limit $b$ such that $ \sigma_n^2/n \to b^2$ then the Heyde's condition (cf\. e.g\. 
\cite{HaHe, Chapter 5}) is satisfied and we get a CLT.
As proved in \cite{VWo14}, in general, for $S_n - E(S_n| \Cal F_0)$ the CLT under Heyde's condition is not quenched. Our theorem shows that it is quenched
in the particular case when $(f\circ T^i)$ is a stationary linear process. 

\underbar{Remark 2}. Theorem 1 implies a quenched CLT for $S_n - E(S_n| \Cal F_0)$ as soon as $\sum_{k=1}^\infty a_k^2 <\infty$, 
$\liminf_{n\to\infty} \bar{\sigma}_n^2/n >0$, and the sequence of $ b_k = \sum_{i=0}^{k} a_i$ is bounded.

\demo{Proof of Theorem 1} 
We have to prove a quenched CLT for the triangular array of random variables $b_{n-k}e\circ T^k/\bar{\sigma}_n$, $k=1\dots,n$, $n=1,2,\dots$.
\smallskip

The $e\circ T^k$ are iid and they remain iid for the conditional probabilities $m_\omega$ as well. From $\bar{\sigma}_n^2= \sum_{k=1}^{n-1} b_k^2
\to \infty$ we get the CLT.
\smallskip

Let $e\circ T^k$ be martingale differences and let (2) hold.
To prove the CLT we use Lachout's refinement \cite{L} of the McLeish's central limit theorem (\cite{Mc}),
applied to regular conditional probabilities with
respect to the $\sigma$-algebra $\Cal F_0$. We thus will prove
\roster
\item"(a)" $E\big(\max_{k\leq n-1} |b_{n-k}e\circ T^k|/\bar{\sigma}_n\, | \,\big|\, \Cal F_0\big)\to 0$ a.s\.,
\item"(b)" $\sum_{k=1}^{n-1} b_{n-k}^2 e^2\circ T^k / \bar{\sigma}_n^2$ converge to a constant a.s..
\endroster

By (2), 
$$
  \frac{b_{n-k}^2}{\bar{\sigma}_n^2} \leq \frac{c}{n}
$$
for all $n$, $1\leq k\leq n-1$, hence (a) follows in the same way as in \cite{VWo14}.

We will prove (b).

Denote 
$$
  T_nf = \frac1{\bar{\sigma}_n^2} \sum_{k=1}^{n-1} b_{n-k}^2 f\circ T^k, \quad f\in L^1.
$$
Recall the Banach's principle (cf\. \cite{K}): \newline
If
\roster
\item"(i)" $T_n$ : $L^1 \to L^1$ are continuous,
\item"(ii)" for every $f\in L^1$, $\sup_n |T_nf| <\infty$ a.e.,
\item"(iii)" there is a dense subset of $h\in L^1$ for which $T_nh$ converges a.s.,
\endroster
then for all $f\in L^1$, $T_nf$ converge a.s.. 
\medskip

We will verify (i)-(iii).
\medskip

(i) follows from the definition.

For (ii),
$$
   |T_nf| \leq \frac1{\bar{\sigma}_n^2} \sum_{k=1}^{n-1} b_{n-k}^2 |f|\circ T^k \leq \frac{c}{n} \sum_{k=1}^{n-1} |f|\circ T^k
$$
hence, by Birkhoff's ergodic theorem, 
$$
  \sup |T_nf| < \infty\,\,\,\,a.s.\quad \forall f\in L^1.
$$
\medskip

Let us prove (iii). Let $f = g-g\circ T$, $g\in L^\infty$. Then
$$\gathered
  T_n f = \frac1{\bar{\sigma}_n^2} \sum_{k=1}^{n-1} b_{n-k}^2 [g\circ T^k - g\circ T^{k+1}] = \\
  = \frac1{\bar{\sigma}_n^2} \sum_{k=1}^n [b_{n-k}^2 - b_{n-k+1}^2] g\circ T^k + \frac{b_n^2}{\bar{\sigma}_n^2} g\circ T \leq \\
  \leq \frac1{\bar{\sigma}_n^2} \sqrt{\sum_{k=1}^n (b_{n-k} + b_{n-k+1})^2} \sqrt{\sum_{k=1}^n a_{n-k+1}^2} \|g\|_\infty 
  + \frac{c}n \|g\|_\infty \leq \\
  \leq \frac2{\bar{\sigma}_n} \sqrt{1+ \frac{c}n} A \|g\|_\infty + \frac{c}n \|g\|_\infty
  \endgathered \tag{3}
$$  
where $A^2 = \sum_{k=1}^\infty a_k^2 <\infty$, hence $T_n f \to 0$ a.s..

The set of functions $c+g-g\circ T$, $c\in \Bbb R$, $g\in L^\infty$, is dense in $L^1$;
for $f' = g-g\circ T$ we have $T_n f \to 0$ a.s. by the calculation above, for $f''=c$ we have $T_nf''=f'' \equiv c$ hence the convergence
towards $c$ takes place for $f=f'+f''$.

By the Banach's principle we conclude that 
$$
  T_ne^2 = \frac1{\bar{\sigma}_n^2} \sum_{k=1}^n b_{n-k}^2 e^2\circ T^k \tag{4}
$$
converges almost surely for every $e\in L^2$. \newline
Let $f^*$ be the limit in (4). Using a similar calculation as in (3) we can see that $T_ne^2 - (T_ne^2)\circ T \to 0$ in $L^1$
hence $f^* = f^*\circ T$. By ergodicity, $f^*$ is a constant a.s..
\enddemo
\qed

\bigskip

\subheading{3. A Non-quenched CLT}

Recall that if the process $(f\circ T^i)$ is adapted and regular and if 
$$
  \sum_{n=1}^\infty \frac{\|E(S_n(f)\,|\,\Cal F_0)\|_2}{n^{3/2}} < \infty,
$$ 
we say that the Maxwell-Woodroofe's condition takes place. If for $P_0f = E(f\,|\, \Cal F_0) - E(f\,|\, \Cal F_{-1})$,
$$
  \sum_{i=0}^\infty \|P_0U^if\|_2 < \infty
$$  
and $(f\circ T^i)$ is adapted and regular then we say that the Hannan's condition takes place.

\proclaim{Theorem 2} There exists a causal stationary linear process $(f\circ T^i)$ with martingale difference innovations such that
\roster
\item"({\bf i})" the Maxwell-Woodroofe's condition and the Hannan's condition
\comment
$$
  \sum_{n=1}^\infty \frac{\|E(S_n(f)\,|\,\Cal F_0)\|_2}{n^{3/2}} < \infty
$$  
\endcomment
are satisfied,
\item"({\bf ii})" for $\sigma_n = \|S_n(f)\|_2$, $\sigma_n \to\infty$, $\sigma_n /\sqrt n \to 0$, $\|E(S_n(f)\,|\,\Cal F_0)\|_2/\sigma_n \to 0$, 
i.e\. $ \bar{\sigma}_n/\sigma_n \to 1$,
%$\|S_n(f) - E(S_n(f)\,|\,\Cal F_0)\|_2/\sigma_n \to 0$,
\item"({\bf iii})" $S_n(f)/\sigma_n$ converge in distribution to the standard normal law $\Cal N(0, 1)$,
\item"({\bf iv})" the convergence is not quenched neither for $S_n(f)/\sigma_n$ nor for $(S_n(f) - E(S_n(f)\,|\,\Cal F_0))/\sigma_n$,
\item"({\bf v})" the WIP does not hold.
\endroster
\endproclaim

\underbar{Remark 3:} When norming by $\sqrt n$, the Hannan's condition implies the WIP.  For $S_n(f) - E(S_n(f)\,|\, \Cal F_0)$ 
the invariance principle is quenched
(cf\. \cite{CuV}), for $S_n(f)$ the CLT is not quenched (cf\. \cite{VWo10}).
\smallskip

\underbar{Remark 4:} As shown in \cite{CuMe}, 
Maxwell-Woodroofe's condition implies a quenched CLT and WIP for $S_n(f)/\sqrt n$ (cf\. also \cite{PU}). 

%In our example, the Maxwell-Woodroofe's condition is 
%satisfied, $S_n(f)/\sqrt n$ converge to zero, $S_n(f)/\sigma_n$ converge to $\Cal N(0,1)$, but the CLT is not quenched and the weak invariance 
%principle does not hold. 

\demo{Proof}
We will find a filtration $(\Cal F_i)_i$ such that $\Cal F_{i+1} = T^{-1}\Cal F_i$ and $e\in L^2(\Cal F_0)\ominus L^2(\Cal F_{-1})$,
$\|e\|_2 =1$. The construction of $e$ and $(\Cal F_i)_i$ will be presented later; it will be needed for the proof of ({\bf iv}) and  
 ({\bf v}) only.

We define a function $f$ by
$$
  f = e + \sum_{k=1}^{\infty} \frac{-\gamma_k}{V_k} \sum_{i=1}^{V_k} U^{-i}e
$$
where $\gamma_k>0$, $\sum_{k=1}^{\infty} \gamma_k = 1$, $V_k \nearrow \infty$, are such that

$$\gather
  \|S_n(f)\|_2\to \infty, \quad \frac{\|S_n(f)\|_2}{\sqrt n} \to 0, \quad
  \frac{\|E(S_n(f)\,|\,\Cal F_0)\|_2}{\|S_n(f)\|_2} \to 0, \\
  \quad\text{and}\quad
  \frac{S_n(f)}{\|S_n(f)\|_2} \to N(0, 1).
  \endgather
$$

To do so, we define
$$
  \gamma_k = \frac2{k+2} \prod_{j=1}^{k} \big(1 - \frac1{j+1}\big), \quad k=1,2,\dots.
$$
From
$$
  \gamma_k = 2 \prod_{j=1}^{k} \big(1 - \frac1{j+1}\big) - 2 \prod_{j=1}^{k+1} \big(1 - \frac1{j+1}\big)
$$
and
$$
  \lim_{k\to\infty} \prod_{j=1}^{k} \big(1 - \frac1{j+1}\big) = 0
$$
we deduce
$$
   \sum_{j=k}^\infty \gamma_j = 2 \prod_{j=1}^{k} \big(1 - \frac1{j+1}\big).
$$
Therefore,
$$
  \sum_{k=1}^\infty \gamma_k =1,\quad 1 - \sum_{j=1}^{k-1} \gamma_j = 2 \prod_{j=1}^{k} \big(1 - \frac1{j+1}\big) = (k+2)\gamma_k.
$$

The numbers $V_k$ will be specified later. We suppose that $V_k$ grow at least exponentially fast.

We have
$$
  \big\|\frac{\gamma_k}{V_k} \sum_{i=1}^{V_k} U^{-i}e \big\|_2 = \frac{\gamma_k}{\sqrt{V_k}}
$$
which guarantees that
$$
  f = \sum_{k=1}^\infty \gamma_k \Big( e - \frac{1}{V_k} \sum_{i=1}^{V_k} U^{-i}e \Big) =
  \sum_{k=1}^\infty \gamma_k f_k \in L^2
$$
where
$$
  f_k = e - \frac{1}{V_k} \sum_{i=1}^{V_k} U^{-i}e =  g_k - Ug_k, \quad g_k = -\frac1{V_k} \sum_{j=1}^{V_k} jU^{-V_k-1+j} e.
$$  

For $h= \sum_{i=0}^\infty c_iU^{-i}e$ we have
$$
  S_n(h) = \sum_{j=0}^{n-1} \sum_{i=0}^\infty c_i U^{j-i}e = \sum_{u=-\infty}^{n-1} \sum_{j=\max\{0, u\}}^{n-1} c_{j-u}U^ue. \tag{5}
$$
Let us denote
$$
  f_k =  e - \frac{1}{V_k} \sum_{i=1}^{V_k} U^{-i}e = \sum_{i=0}^\infty c_{k,i}U^{-i}e.
$$
We then have
$$\gathered
  \big| \sum_{j=0\vee u}^{n-1} c_{k,j-u}\big| \leq 1\,\,\,\,\text{for every}\,\,\,\, u, \\
  \big| \sum_{j=0\vee u}^{n-1} c_{k,j-u}\big| \leq \frac{n}{V_k} \,\,\,\,\text{for}\,\,\,\,-1\geq u \geq -V_k, \\
  \big| \sum_{j=0\vee u}^{n-1} c_{k,j-u}\big| = 0 \,\,\,\,\text{for}\,\,\,\,u<-V_k \\
  \sum_{j=0\vee u}^{n-1} c_{k,j-u} \geq 0 \vee 1 - \frac{n}{V_k} \,\,\,\,\text{for}\,\,\,\, u\geq 0.
  \endgathered \tag6
$$  
We deduce that for every $V_k$,
$$
  \|S_n(f_k)\|_2 = \big\|S_n\big( e - \frac{1}{V_k} \sum_{i=1}^{V_k} U^{-i}e \big)\big\|_2 \leq \sqrt{2 n}
$$
hence by the Lebesgue Dominated Convergence Theorem
$$
  \frac{\|S_n(f)\|_2}{\sqrt n} \leq \sum_{k=1}^\infty \gamma_k \frac{\|S_n(f_k)\|_2}{\sqrt n} \to 0. \tag7
$$

\noindent Recall that
$$
  1 - \sum_{j=1}^{k+1} \gamma_j = \sum_{j=k+2}^\infty \gamma_j = (k+4) \gamma_{k+2}.
$$
For $V_k \leq n <V_{k+1}$ we using (5), (6), and the (at least) exponential growth of the $V_k$ get
$$\multline
  \|S_n(f)\|_2 \geq \|S_n(f) - E(S_n(f)\,|\,\Cal F_0)\|_2 \geq \sqrt{n-1}  \big(1 - \sum_{j=1}^{k+1} \gamma_j -
  \sum_{j=k+2}^{\infty} \frac{V_{k+1}}{V_j}\gamma_j \big) \geq \\
  \geq C \sqrt n  (k+4) \gamma_{k+2} 
  \endmultline \tag8
$$
for some constant $C>0$.
Supposing
$$
  \sqrt{V_k} \big(1 - \sum_{j=1}^{k+1} \gamma_j \big) \to \infty
$$
we thus get
$$
  \|S_n(f)||_2 \to \infty. \tag9
$$
Using the (at least) exponential growth of the $V_k$ again we have, for $V_k \leq n <V_{k+1}$,
$$
  \|E(S_n(f)\,|\,\Cal F_0)\|_2 \sim (\gamma_k+\gamma_{k+1})\sqrt n. \tag{10}
$$
To prove this, first recall that $f_j = g_j-Ug_j$ where $g_j\in L^2$. If $V_k$ is large enough we thus get 
$\|E(S_n(\sum_{j=1}^{k-1} f_j)\,|\,\Cal F_0)\|_2 =o(n)$. For $V_{k+2}$ large enough we get
$\|E(S_n(\sum_{j=k+2}^\infty f_j)\,|\,\Cal F_0)\|_2 =o(n)$ by (6).\newline
Using (5) and (6) we can see that
$$
  \|E(S_n(f_k) \,|\,\Cal F_0)\|_2^2 = \frac1{V_k^2} \sum_{j=1}^{V_k} j^2 \sim V_k
$$
and
$$
  \|E(S_n(f_{k+1}) \,|\,\Cal F_0)\|_2^2 \leq  V_{k+1}+1.
$$
Because $\gamma_k \sim \gamma_{k+1}$ we get (10).

By definition,
$$
  \gamma_k = \frac1{k+2} \big(1 - \sum_{j=1}^{k-1} \gamma_j \big)
$$
hence by (8) and (10) we have
$$
  \frac{\|E(S_n(f)\,|\,\Cal F_0)\|_2}{\|S_n(f)\|_2} \to 0. \tag{11}
$$
From (9), (7), and (11) we get ({\bf ii}).
\medskip

By \cite{PU06} and $\sigma_n^2 \to \infty$ we get the central limit theorem ({\bf iii}).
\bigskip

\centerline{\it The Hannan's and Maxwell-Woodroofe's condition}
\medskip

We have $\|P_0 U^if\|_2 = |a_i|$, $i\geq 0$, where $f = a_0e -\sum_{i=1}^\infty a_iU^{-i}e$. From the definition of $f$
we deduce that $a_0=1$ and $a_i>0$ for $i\geq 1$, $\sum_{i=1}^\infty a_i = 1$. This implies the Hannan's condition.
\medskip

Denote 
$$
  h_k = \frac{-\gamma_k}{V_k} \sum_{i=1}^{V_k} U^{-i}e, \quad f' = \sum_{k=1}^\infty h_k;
$$  
we thus have $f=e-f'$. By (6), for any $k\geq 1$, 
$$
  \|S_n(h_k)\|_2^2 \leq \gamma_k^2 2V_k \Big(\frac{n}{V_k}\Big)^2 = 2 \gamma_k^2 \frac{n^2}{V_k}, \quad n=1,2,...,V_k
$$
and
$$
  \|E(S_n(h_k)\,|\,\Cal F_0)\|_2 = \|E(S_{V_k}(h_k)\,|\,\Cal F_0)\|_2 \leq \gamma_k \sqrt{V_k} , \quad n\geq V_k
$$
hence
$$
  \sum_{n=1}^\infty \frac{\|E(S_n(h_k)\,|\,\Cal F_0)\|_2}{n^{3/2}} \leq \sqrt 2 \gamma_k \frac{1}{\sqrt{V_k}} 
  \sum_{n=1}^{V_k} \frac{1}{\sqrt{n}} +
  \gamma_k \sqrt{V_k} \sum_{n=V_k+1}^\infty \frac1{n^{3/2}} \leq C\gamma_k
$$
for some constant $C$. Therefore,
$$
  \sum_{n=1}^\infty \frac{\|E(S_n(f')\,|\,\Cal F_0)\|_2}{n^{3/2}} \leq
  \sum_{k=1}^\infty \sum_{n=1}^\infty \frac{\|E(S_n(h_k)\,|\,\Cal F_0)\|_2}{n^{3/2}}
  \leq C \sum_{k=1}^\infty \gamma_k = C <\infty
$$
and ({\bf i}) follows.

\bigskip

\centerline{\it The filtration and $e$}
\medskip

For all what have been proved up to now we supposed only that $(e\circ T^i)_i$ are martingale differences and that 
$\|e\|_2 =1$. In order to get ({\bf iv}) and ({\bf v}) we will need a particular choice of the filtration and of $e$.

Let $\Cal B_k',\Cal B_l'' \subset \Cal A$, $k, l=1,2,\dots$, be mutually independent $\sigma$-algebras,
$$
%\gather
  \Cal B_k'\subset T^{-1}\Cal B_k',\quad \cap_{j=1}^\infty T^{j}\Cal B_k' = \{\Omega, \emptyset\},\quad 
  \Cal B_l'' =  T^{-1}\Cal B_l''
%',\\
 % \cap_{j=1}^\infty T^{j}\Cal B_k' = \{\Omega, \emptyset\},\quad  \cap_{j=1}^\infty T^{j}\Cal B_l'' = \{\Omega, \emptyset\}
 % \endgather
$$
(modulo sets of measure 0 or 1) for every $k, l$;
$\xi_k\circ T^i$ are iid $\Cal B_k'$-measurable random variables, $\mu(\xi_k=1) = 1/2 = \mu(\xi_k=-1)$ for all $i$.

All these objects can be constructed by taking finite alphabets $\Bbb A_k'$ and $\Bbb A_l''$, $k,l=1,2,\dots$,
$\Omega_k' = \underset i\in \Bbb Z \to{\times} \Bbb A_{k,i}'$ where $\Bbb A_{k,i}'$ are identical copies of $\Bbb A_k'$, similarly we define
$\Omega_l''$, $k,l=1,2,\dots$. On the sets $\Omega_k'$ and $\Omega_l''$ we define product $\sigma$-algebras, product measures, and 
left shift transformations $T_k'$, $T_l''$. $\Omega$ is the product of all $\Omega_k'$ and $\Omega_l''$ equipped with the product 
$\sigma$-algebra $\Cal A$, the product (probability) measure $\mu$, and the product transformation $T$. 
For projections $\xi_k$ and $\zeta_l$ of $\Omega$ onto $\Bbb A_{k,0}'$ and $\Bbb A_{l,0}''$ we thus get mutually independent processes of 
iid $(\xi_k\circ T^i)_i$, $(\zeta_l\circ T^i)_i$. We suppose that $\Bbb A_k' = \Bbb A_k'' = \{-1, 1\}$, $k=1,2, \dots$, and 
$\mu(\xi_k=1) = 1/2 = \mu(\xi_k=-1) = \mu(\zeta_k=1) = \mu(\zeta_k=-1)$. For $\Cal B_k'$ we take the past $\sigma$-algebras
$\sigma\{\xi_k\circ T^i : i\leq 0\}$ and for $\Cal B_l''$ we take the $\sigma$-algebras $\sigma\{\zeta_l\circ T^i : i\in \Bbb Z\}$.
The properties above can be easily verified, the latter follow from Kolmogorov's 0-1 law.

We thus have 
that $\xi_k\circ T^i$ are iid $T^{-i}\Cal B_k'$-measurable random variables, $\mu(\xi_k=1) = 1/2 = \mu(\xi_k=-1)$.

Recall that $ \sigma_n/\sqrt n = \|S_n(f)\|_2 /\sqrt n \to 0$ (cf\. (7)).
We thus can suppose that $N_k$ is big enough so that for $k$ odd, 
$$ 
   2^k\sigma_{N_k} \leq  \frac{\sqrt{N_k}}{k^{3/2}}\,\,\,\,\text{ and}\,\,\,\, N_{k+1} = 4N_k, \,\,\,\,\sigma_{4N_k} \leq 
  2 \sigma_{N_k}, \,\,\,\, \sum_{k=1}^\infty \frac1{4N_k} < \frac12. \tag{12}
$$

For $k=1,2,\dots$, let $A_k\in \Cal B_k''$ be sets such that $T^{-i}A_k$, $i=0,\dots,3N_k$ are mutually disjoint 
(hence $\{T^{-i}A_k :\, i=0,\dots,3N_k\}$ are Rokhlin towers) and $\mu(A_k)= 1/(4N_k)$
%the values of $N_k$ will be specified later 
(existence of Rokhlin Towers is proved e.g\. in \cite{CSF}). From (12) if follows
$$
  \sum_{k=1}^\infty \mu(A_k) < \frac12. \tag{13}
$$

By $\Cal B''$ we define the $\sigma$-algebra generated by all $\Cal B_k''$; we thus have $T^{-1}\Cal B''
=\Cal B''$, all Rokhlin towers defined above are $\Cal B''$-measurable.\newline
By $\Cal F_j$ we denote the $\sigma$-algebra generated by $\Cal B''$ and all $\xi_k\circ T^i$, $i\leq j$, $k=1,2,\dots$; notice that
$T^{-1}\Cal F_j = \Cal F_{j+1}$.

For $d = 2 (\sum_{k=1}^\infty 1/{k^3})^{-1/2}$ we define
$$
  e_k = d\xi_k\frac{\sqrt{N_k}}{k^{3/2}}1_{A_k}, \quad e = \sum_{k=1}^\infty e_k.
$$
Notice that $\|e_k\|_2 = d/(2 k^{3/2})$ hence $e\in L^2$. By definition, $e$ is $\Cal F_0$-measurable. 
By definition, $e_k$ are mutually independent hence $\|e\|_2^2 = (d^2/4) \sum_{k=1}^\infty 1/{k^3}$; we thus have 
$\|e\|_2 = 1$. \newline
Because $A_k\in \Cal B''$ and $\xi_k$ is independent of $\Cal F_{-1}$, we have $E(e_k\,|\,\Cal F_{-1}) =0$ for every $k$ hence 
$E(e\,|\,\Cal F_{-1}) =0$, $(U^ie)_i$ is thus a martingale difference sequence adapted to the filtration $(\Cal F_i)$.

Recall that
$$
  f = e + \sum_{k=1}^{\infty} \frac{-\gamma_k}{V_k} \sum_{i=1}^{V_k} U^{-i}e = a_0e - \sum_{i=1}^\infty a_i U^{-i}e
$$
where $a_0=1$, $a_i>0$ for all $i\geq 1$, and $\sum_{i=1}^\infty a_i = 1$.

By $m_{\omega}$ we will denote regular conditional probabilities w.r.t\. $\Cal F_0$ ($\Cal A$ is a Borel $\sigma$-algebra of a Polish space
hence the regular conditional probabilities exist). Notice that all sets $T^{-i}A_k$, $k=1,2,\dots$, 
$i\in\Bbb Z$, are $\Cal F_0$-measurable hence $m_{\omega}(T^{-i}A_k) = 0$ (if $\omega\not\in T^{-i}A_k$) or $m_{\omega}(T^{-i}A_k) = 1$ 
(if $\omega\in T^{-i}A_k$). 
\smallskip

Let us fix a $k\geq 1$ odd 
%such that
%$$
%  \sum_{i=N_k}^\infty a_i < \frac12 
%$$
and denote 
$$
  A_k' = A_k \setminus \underset j {\neq k} \to{\bigcup} A_j.
$$
By (13) and independence, $\mu(A_k') \geq \mu(A_k)/2$.
The sets $A_k',\dots, T^{-3N_k+1}A_k'$ are mutually disjoint and $\mu(A_k)\geq 1/(4N_k)$ hence
$$
  \mu\Big(\underset N=0 \to{\overset N_k-1\to{\bigcup}} T^{-N+1} A_k'\Big) \geq \frac1{8}.\tag{14}
%  \big(A_k \setminus \underset j\neq k \to{\cup} A_j\big) \Big) \geq \frac1{16}.\tag{13}
$$  

We have 
$$\gathered
  S_N(f) = \sum_{j=0}^{N-1} U^j\Big(e - \sum_{i=1}^\infty a_iU^{-i}e\Big) = \\
   U^{N-1}e + \sum_{j=1}^{N-2} U^je - \sum_{i=0}^{\infty}  \sum_{j=0\vee 1-i}^{N-1} a_{i+j}U^{-i}e + 
   e - \sum_{i=2-N}^{-1}  \sum_{j=1-i}^{N-1} a_{i+j}U^{-i}e = \\
%  U^{N-1}e + \sum_{j=0}^{N-2} U^je_k -  \sum_{j=0}^{N-1}  \sum_{i=1}^\infty a_iU^{j-i}e_k +
% \sum_{l\geq 1, l\neq k}\sum_{j=0}^{N-2} U^je_l -  
%\sum_{l\geq 1, l\neq k}\sum_{j=0}^{N-1}  \sum_{i=1}^\infty a_iU^{j-i}e_l =\\
  U^{N-1}e + I -II +III -IV.
  \endgathered \tag{15}
$$

Let us suppose that  $N_k\leq N < N_{k+1} = 4N_k$ ($k$ is odd).
For $\omega \in T^{-N+1} A_k'$
%\Big(A_k \setminus \underset j \neq k \to{\bigcup} A_j\Big)$
we have
%$$
%  |U^{N-1}e(\omega)| = |U^{N-1}e_k(\omega)| = d\frac{\sqrt{N_k}}{k^{3/2}} 1_{T^{-N+1}A_k}(\omega) = 
% d\frac{\sqrt{N_k}}{k^{3/2}}
%$$
%and
$$
  m_\omega\big(U^{N-1}e = d\frac{\sqrt{N_k}}{k^{3/2}}\big) = m_\omega\big(U^{N-1}e = -d\frac{\sqrt{N_k}}{k^{3/2}}\big) 
  = \frac12.
$$
$-II+III = E(S_N(f) | \Cal F_0)$ hence it is a constant $m_\omega$ almost surely. \newline
$U^{N-1}e + I-IV$ is (an infinite) linear combination of products of $U^i\xi_l$  with $\Cal F_0$-measurable functions, 
$1\leq i\leq N-1$.
Because $\Cal F_0$,  $U^i\xi_l$, $1\leq i\leq N-1$, are mutually independent,  $U^i\xi_l$, $1\leq i\leq N-1$, are iid 
with respect to the measure $m_\omega$ and $m_\omega(\xi_l =\pm 1) = 1/2$ (for $\mu$ a.e\. $\omega$).
Therefore, $I-IV$ is a symmetric random variable independent of $U^{N-1}e$ w.r.t\. $m_\omega$; $U^{N-1}e + I-IV
= S_N(f) -  E(S_N(f) | \Cal F_0)$ is a symmetric random variable as well.

\comment
$III -IV$ is independent of $U^{N-1}e_k$ and it is a linear combination of mutually independent and symmetric random variables.
Therefore, 
$$
   m_\omega\big( |U^{N-1}e +I-II| \geq  d\frac{\sqrt{N_k}}{k^{3/2}}\big) \geq \frac14.
$$

For $l=k$, after $1_{T^{-u}A_l}=1$ we have $1_{T^{-u}A_l}=0$ for at least $3N_k$ consecutive numbers $u$; we thus 
on $ T^{-N+1} A_k'$ have
$I =0$ and
$$\gather
   1_{T^{-N+1}A_k} \Big|  \sum_{j=0}^{N-1}  \sum_{i=1}^\infty a_iU^{j-i}e_k \Big| = \\
    1_{T^{-N+1}A_k} \Big| \sum_{i=2-N}^\infty \sum_{j=1-i \vee 0}^{N-1} a_{j+i} U^{-i} e_l   \Big| \leq
    1_{T^{-N+1}A_k} |e_k| \sum_{i=N_k}^\infty a_i < \frac12  1_{T^{-N+1}A_k} |e_k|.
   \endgather
$$
\endcomment
We thus have 
$$
    m_\omega\big( |S_N(f) -  E(S_N(f) | \Cal F_0)|  \geq  d\frac{\sqrt{N_k}}{k^{3/2}}\big) =   
   m_\omega\big( |U^{N-1}e +I -IV| \geq 
    d\frac{\sqrt{N_k}}{k^{3/2}}\big) \geq \frac12.
$$
%Because  $N_k\leq N < N_{k+1} = 4N_k$ 
From (15) and $\sum_{i=1}^\infty a_i = 1$ we by direct calculation deduce that $\sigma_n \leq \sigma_{n+1}$ for all $n\geq 1$.
Using (12) we thus get
$$\multline
    m_\omega\Big( \frac{|S_N(f) -  E(S_N(f) | \Cal F_0)|}{\sigma_N}  \geq 2^{k}\frac{\sigma_{N_k}}{\sigma_{N}}\Big) \geq \\
     m_\omega\Big( \frac{|S_N(f) -  E(S_N(f) | \Cal F_0)|}{\sigma_N}  \geq 2^{k-1}\Big) \geq \frac12
     \endmultline
$$
Because $ E(S_N(f) | \Cal F_0)$ is $m_\omega$ a.s\. a constant and $ S_N(f) -  E(S_N(f) | \Cal F_0)$ is a symmetric random 
variable, we get 
$$
    m_\omega\Big( \frac{|S_N(f) | \Cal F_0)|}{\sigma_N}  \geq 2^{k-1}\Big) \geq \frac14.
$$
\comment
and on $T^{-N+1} A_k'$,
$$\gathered
%  S_N(f) = \sum_{j=0}^{N-1} \sum_{i=0}^\infty a'_i U^{j-i}e = \sum_{u=-\infty}^{N-1} \sum_{j=\max\{0, u\}}^{N-1} a'_{j-u}U^ue
  S_N(f) = \sum_{j=0}^{N-1} \sum_{i=0}^\infty a_i U^{j-i}e = \sum_{u=-\infty}^{N-1} \sum_{j=0\vee u}^{N-1} a_{j-u}U^ue =   \\
  = U^{N-1}e_k
  - d \frac{\sqrt{N_k}}{k^{3/2}} \sum_{u=-\infty}^{0} \sum_{j=0\vee u}^{N-1} a_{j-u} 1_{\{T^{-u}A_k\}} U^u\xi_k - \\
  - d \frac{\sqrt{N_k}}{k^{3/2}} \sum_{u=1}^{N-2} \sum_{j=0\vee u}^{N-1} a_{j-u} 1_{\{T^{-u}A_k\}} U^u\xi_k - \\
  - d \sum_{l\neq k} \frac{\sqrt{N_l}}{l^{3/2}} \sum_{u=-\infty}^{0} \sum_{j=0\vee u}^{N-1} a_{j-u}  1_{\{T^{-u}A_l\}}U^u\xi_l -\\
  - d \sum_{l\neq k} \frac{\sqrt{N_l}}{l^{3/2}} \sum_{u=1}^{N-2} \sum_{j=0\vee u}^{N-1} a_{j-u}  1_{\{T^{-u}A_l\}}U^u\xi_l = \\
  = U^{N-1}e_k - I - II - III-IV.
  \endgathered \tag{14}
$$
On $T^{-N+1} A_k'$, $\omega\in T^{-N+1} A_k'$, we have 
$$
  m_\omega\big(U^{N-1}e = d\frac{\sqrt{N_k}}{k^{3/2}}\big) = 
   m_\omega\big(U^{N-1}e = -d\frac{\sqrt{N_k}}{k^{3/2}}\big) = \frac12.
$$
$U^{N-1}e_k$ is independent of $I$, $II$, $III$, $IV$ (with respect to $m_\omega$, $\omega \in T^{-N+1} A_k'$) hence
by symmetry of $I+II+III+IV$ the probability $m_\omega$ that $U^{N-1}e_k$ and $I+II+III+IV$ have the same sign is at least $1/4$,
therefore
$$
  m_\omega\big( |S_N(f)| \geq d\frac{\sqrt{N_k}}{k^{3/2}}\big) \geq \frac14. \tag{15}
$$  
Recall that 
$$
  \sigma_n/\sqrt n = \|S_n(f)\|_2 /\sqrt n \to 0.
$$
We thus can suppose that $N_k$ is big enough so that  $2^k\sigma_{N_k} \leq d {\sqrt{N_k}}/{k^{3/2}}$; at the same time we can suppose
that $k$ is odd and that $N_{k+1} = 4N_k$.
From this and from (15) it follows that there exists a $\delta>0$ such that for $k$ sufficiently big and $N_{k+1} = 4N_k$, 

We deduce that for $k$ odd and bigger than some $k_0$ there exists a $\delta>0$ not depending on $k$ and $N$ such that
the Prokhorov (metric) distance  
between the distribution of $(S_N(f)- E(S_N(f) | \Cal F_0)) /\sigma_N$ and $\Cal N(0, 1)$ is bigger than $\delta$.
\endcomment
% For the given odd $k$ there thus 
For any $K<\infty$ there thus exists a $k_0$ (with $2^{k_0-1} \geq K$) such that for an integer $k\geq k_0$ there
exists a set $B_k = \cup_{N_k}^{N_{k+1}-1} T^{-N+1} A_k'$ of measure bigger than $1/16$ (cf\. (14)) such 
that for $\omega\in B_k$ and the probability $m_\omega$ there exists
an $N_k\leq N< N_{k+1}$ for which 
%the Prokhorov (metric) distance  between the distribution of 
%$(S_N(f)- E(S_N(f) | \Cal F_0))/\sigma_N$ (i.e\. 
%$m_\omega\circ ( (S_N(f)- E(S_N(f) | \Cal F_0) )/\sigma_N)^{-1}$) and $\Cal N(0, 1)$ is bigger than $\delta$.
$$
   m_\omega\Big( \frac{|S_N(f) -  E(S_N(f) | \Cal F_0)|}{\sigma_N}  \geq K\Big) \geq \frac12, \quad
   m_\omega\Big( \frac{|S_N(f) | \Cal F_0)|}{\sigma_N}  K\Big) \geq \frac14. \tag{16}
$$
We conclude that there exists a set $B$ of positive measure such that for $\omega\in B$ there is an infinite sequence of
$k$ (odd) and $N_k\leq N\leq N_{k+1}$ such that for the probability $m_\omega$, the laws of 
$(S_N(f) - E(S_N(f) | \Cal F_0))/\sigma_N$ and of $S_N(f)/\sigma_N$ do not weakly converge to
$\Cal N(0, 1)$. This proves that the CLT for $(S_n - E(S_n | \Cal F_0))/\sigma_n$ and for $S_N(f)/\sigma_N$ are not quenched. 

\comment
Let us denote $X_n = (S_n - E(S_n | \Cal F_0))/\sigma_n$, $a_n =   E(S_n | \Cal F_0)/\sigma_n$. For ($\mu$ a.e\.) probability 
$m_\omega$, $a_n$ is ($m_\omega$) a.s\. a constant and by the considerations following formula (15), $X_n$ is a symmetric
random variable (note that $U^{N-1}e$ is independent of $I-IV$); we speak of properties w.r.t\. the probability 
$m_\omega$. \newline
If the distributions $m_\omega\circ X_n^{-1}$ do not weakly converge to $\Cal N(0, 1)$, $m_\omega\circ (X_n+a_n)^{-1}$ 
do not weakly converge to $\Cal N(0, 1)$ either.  \newline
To see this, denote $\varphi_n(t) = E(\exp(it(X_n+a_n)))$ the characteristic
function of $X_n+a_n$. If $m_\omega\circ (X_n+a_n)^{-1}$ weakly converge to $\Cal N(0, 1)$ then 

$\varphi_n(t) = E(\exp(it(X_n)) \exp({ita_n}) \to e^{-t^2/2}$ \newline 
for every $t\in\Bbb R$. Because $X_n$ are symmetric,
$ E(\exp(it(X_n))$ is an even function of $t$ hence for every $t\in\Bbb R\setminus \{0\}$, $ \exp({ita_n}) \to \pm 1$.
Therefore, $a_n\to 0$ hence $m_\omega\circ X_n^{-1}$ weakly converge to $\Cal N(0, 1)$, in contradiction to our assumptions.

The CLT for $S_n/\sigma_n$ is thus not quenched. 
\endcomment
This finishes the proof of  ({\bf iv}).

 \comment

From (15) and the the text 
after it follows that   From this it follows that for ($\mu$ a.e\.) probability $m_\omega$, $S_n - E(S_n | \Cal F_0)$ is a symmetric
random variable. Therefore, if $m_\omega\circ ((S_n - E(S_n | \Cal F_0))/\sigma_n)^{-1}$ do not converge to $\Cal N(0, 1)$ then
$m_\omega\circ (S_n/\sigma_n)^{-1}$ do not converge to $\Cal N(0, 1)$ either hence
the CLT for $S_n/\sigma_n$ is not quenched. We thus have proved  ({\bf iv}).
%\comment
\smallskip

In (14) we can note that $E(S_N | \Cal F_0) = -I-III$ and we deduce in the same way as above that the CLT for $S_n - E(S_n | \Cal F_0)$
is not quenched either. This proves ({\bf iv}).
\endcomment
\bigskip

\comment
We have proved that there exists a $C>0$ such that 
$$
  \mu(\max_{1\leq n\leq N_{k+1}} |S_n(f)| > C2^{k} \sigma_{N_{k+1}}) > \frac1{32}
$$
for infinitely many $k$ hence the Donsker invariance principle is not satisfied (more precisely, by \cite{PU05} it is satisfied 
for $S_n(f)/\sqrt n$ with a convergence towards a degenerate limit, but it is not satisfied for $S_n(f)/\sigma_n$).
\endcomment

We have proved (cf\. (16)) that for any $K<\infty$, $k$ sufficiently big (with $2^{k-2} \geq K$), and $\omega\in B_k = 
\cup_{N_k}^{N_{k+1}-1} T^{-N+1} A_k'$ ($B_k$ of measure bigger than $1/16$),
$$
   m_\omega \Big( \underset N=N_k \to{\overset N_{k+1}-1 \to{\cup}} \Big\{   \frac{|S_N(f) -  E(S_N(f) | \Cal F_0)|}
  {\sigma_{N_{k+1}}}  \geq K  \Big \} \Big) \geq \frac12
$$
hence 
$$\multline
   \mu \Big( \underset N=N_k \to{\overset N_{k+1}-1 \to{\cup}} \Big\{   \frac{|S_N(f) -  E(S_N(f) | \Cal F_0)|}
  {\sigma_{N_{k+1}}}  \geq K  \Big \} \Big) =\\
    \mu \Big(  \max_{N_k\leq N\leq N_{k+1}-1}  \Big\{   \frac{|S_N(f) -  E(S_N(f) | \Cal F_0)|}
  {\sigma_{N_{k+1}}}  \geq K  \Big \}  \geq \frac1{32},
   \endmultline
$$
similarly we get
$$
    \mu \Big(  \max_{N_k\leq N\leq N_{k+1}-1} \Big\{   \frac{|S_N(f) -  E(S_N(f) | \Cal F_0)|}
  {\sigma_{N_{k+1}}}  \geq K  \Big \} \Big) \geq \frac1{64}.
$$
This  proves ({\bf v}).

\enddemo
\qed

\bigskip

\underbar{Acknowledgement.} The authors thank Professor Michael Lin for many helpful consultations.

\enddocument
\end

\ref \key B \by Billingsley, P. \paper The Lindeberg-L\'evy theorem for martingales \jour Proc. Amer. Math. Soc.
\vol 12 \pages 788-792 \yr 1961 \endref 
\ref \key C-F-S \by Cornfeld, I.P., Fomin, S.V., and Sinai, Ya.G. \book Ergodic Theory \publ Springer \publaddr Berlin \yr 1982 \endref
\ref \key Cu 1 \by Cuny, Ch. \paper Pointwise ergodic theorems with rate and application to limit theorems for stationary processes
\paperinfo sumitted for publication, arXiv:0904.0185 \yr 2009 \endref
\ref \key Cu 2 \by Cuny, Ch. \paper Norm convergence of some power-series of operators in $L^p$ with applications in ergodic theory
\paperinfo sumitted for publication \yr 2009 \endref
\ref \key D-L 1 \by Derriennic, Y. and Lin, M. \paper Sur le th\'eor\`eme limite central de Kipnis et Varadhan pour les cha\^\i nes 
r\'eversibles ou normales \jour CRAS 323 \pages 1053-1057 \yr 1996 \endref 
\ref \key D-L 2 \by Derriennic, Y and Lin, M. \paper The central limit theorem for Markov
chains with normal transition operators, started at a point \jour Probab. Theory Relat. Fields
\vol 119 \pages 509-528 \yr 2001 \endref
\ref \key G \by Gordin, M.I. \paper A central limit theorem for stationary processes \jour Soviet Math. Dokl. \vol 10 
\pages 1174-1176 \yr 1969 \endref
\ref \key G-Ho \by Gordin, M.I. and Holzmann, H. \paper The central limit theorem for stationary Markov chains under invariant splittings 
\jour Stochastics and Dynamics \vol 4 \pages 15-30 \yr 2004 \endref 
\ref \key G-L 1 \by Gordin, M.I. and Lif\v sic, B.A. \paper Central limit theorem for statioanry processes \jour Soviet Math. Doklady 19 
\pages 392-394 \yr 1978 \endref 
\ref \key G-L 2 \by Gordin, M.I. and Lif\v sic, B.A. \paper A remark about a Markov process with normal transition operator
\paperinfo In: Third Vilnius Conference on Probability and Statistics \vol 1 \pages 147-148 \yr 1981 \endref
\ref \key G-L 3 \by Gordin, M.I. and Lif\v sic, B.A. \paper Th\input amstex
\documentstyle{amsppt}
\magnification 1200

\NoBlackBoxes
\NoRunningHeads
\topmatter
\title Quenched central limit theorems for a stationary linear process
\endtitle
\author Dalibor Voln\'y and Michael Woodroofe\endauthor
\abstract
We find a sufficient condition under which a central limit theorem for a a stationary linear process is quenched.
We find a stationary linear process satisfying the Maxwell-Woodroofe's condition for which $\sigma_n=\|S_n\|_2 = o(\sqrt n)$,
$S_n/\sigma_n$ converge to the standard normal law, and the convergence is not quenched.
\endabstract
\endtopmatter
\document

\subheading{1. Introduction} 

Let $T$ be an ergodic automorphism of a probability space $(\Omega, \Cal A, \mu)$.  
For $h\in L^2$, $Uh= h\circ T$ is a unitary operator; we will freely switch from the notation $h\circ T^i$ to $U^ih$
and vice versa.

Let $(\Cal F_i)_i$ be
a filtration such that $\Cal F_{i+1} = T^{-1}\Cal F_i$, and $e\in L^2(\Cal F_0)\ominus L^2(\Cal F_{-1})$.
For simplicity we will suppose $\|e\|_2 = 1$.
Let $a_i$ be real numbers with $\sum_{i\in \Bbb N} a_i^2 <\infty$ and let
$$
  f = \sum_{i\in \Bbb Z} a_{-i} U^ie.
$$
Then $f\in L^2$ and we say that $(f\circ T^i)_i$ is a stationary linear process. 
The stationary linear process is a classical and important case of a (strictly) stationary process and, moreover, any regular stationary 
process is a sum of stationary linear process ``living'' in mutually orthogonal and $U$-invariant subspaeces of $L^2$: 
If $e_k\in  L^2(\Cal F_0)\ominus L^2(\Cal F_{-1})$, $\|e_k\|_2 = 1$ are  mutually orthogonal, 
$a_{k,i}$ are real numbers with $\sum_{k=1}^\infty \sum_{i\in \Bbb N} a_{k,i}^2 <\infty$, and if
$$
  f = \sum_{k=1}^\infty \sum_{i\in \Bbb N} a_{k,-i} U^ie_k \tag{1}
$$
then we say that $(f\circ T^i)_i$ is a superlinear process. As shown in \cite{VWoZ}, if $f\in L^2$ is $\Cal F_0$-measurable
and $E(f | \Cal F_{-\infty}) = 0$ then a representation (1) exists. Any adapted ($f$ is $\Cal F_0$-measurable) and regular
($E(f | \Cal F_{-\infty}) = 0$) $L^2$ stationary process $(f\circ T^i)_i$ is thus a sum of mutually orthogonal stationary linear processes.
%Stationary linear processesare thus building blocs of general regular processes.

Let us denote $S_n(f) = \sum_{i=0}^{n-1} f\circ T^i$.
Recall (\cite{PU06}) that if $\sigma_n = \|S_n(f)\|_2 \to \infty$ then the distributions of $S_n(f)/\sigma_n$ weakly
converge to $N(0,1)$, i.e\. we have a CLT.

We will study when this CLT is quenched.

Let us suppose that the regular conditional probabilities $m_\omega$ with respect to the $\sigma$-field $\Cal F_0$ exist. If for 
$\mu$ a.e\. $\omega$ the distributions $m_\omega(S_n(f)/\sigma_n)^{-1}$ weakly converge to $N(0,1)$, we say that the CLT is quenched.
A quenched CLT can be defined using Markov Chains. Any stationary process $(f\circ T^i)_i$ can be expressed using a homogeneous and 
stationary Markov Chain $(\xi_i)_i$ as $(g(\xi_i))_i$; a CLT is quenched if it takes place for a.e\. starting point (this approach 
is probably earlier than our one; it has been used in e.g\. \cite{DLi}). Equivalence of both approaches was explained e.g\. 
in \cite{V}.

In the next section, for a stationary linear processes, we will give a sufficient condition for a quenched CLT. Then, in Section 2,
we will present a stationary linear processes $(f\circ T^i)_i$ for which the CLT is not quenched. As shown in \cite{CuP}, 
Maxwell-Woodroofe's condition implies a quenched CLT for $S_n(f)/\sqrt n$. In our example, the Maxwell-Woodroofe's condition is satisfied, 
but the CLT is not quenched for $S_n(f)/\sigma_n$. The weak invariance principle is not satisfied either. \newline
An example of a non quenched CLT has been found
in \cite{VWo10} for a stationary linear process (*) satisfying Hannan's condition (cf\. \cite{HaHe}), i.e\. $
\sum_{i\in \Bbb N} |a_{i}| <\infty$. Under Hannan's condition, the CLT is quenched for $S_n(f) - E(S_n(f) | \Cal F_0)$, cf\.
\cite{CuP} and \cite{VWo}. In our example, the CLT is not quenched neither for $S_n(f) - E(S_n(f) | \Cal F_0)$ nor for
$S_n(f)$.
\bigskip

\subheading{1. A Sufficient Condition}
Let $(e\circ T^i)$ be a martingale difference sequence as defined in the introduction,
%$e_k = e\circ T^k$ be a sequence of martingale differences adapted to a filtration $(\Cal F_i)_i$, $\Cal F_i = T^{-i}\Cal F_0$, 
$$
  f = \sum_{i=0}^\infty a_ie\circ T^{-i}
$$
where $\sum_{i=0}^\infty  a_i^2 < \infty$.
%$(f\circ T^k)$ is thus a stationary linear process. We denote $e_k = e\circ T^k$ 
and $S_n = \sum_{i=0}^{n-1} f\circ T^i$, $n=1,2,\dots$. By definition,
$$
  S_n = \sum_{j=0}^{n-1} f\circ T^j = \sum_{i=0}^\infty \sum_{j=0}^{n-1} a_ie\circ T^{j-i} = \sum_{k=-\infty}^{n-1} \sum_{j=k}^{n-1} 
  a_{j-k}e\circ T^k
  = \sum_{k=-\infty}^{n-1} b_{n-k} e\circ T^k.
$$  
where
$$
  b_0=0,\,\,\, b_j = \sum_{i=0}^{j-1} a_i,\quad j\geq 1.
$$  

We denote 
$$
  \sigma_n^2 = E[(S_n - E(S_n| \Cal F_0))^2] = \sum_{k=1}^{n-1} b_{n-k}^2.
$$  

\proclaim{Theorem 1} Let $\sigma_n^2 \to \infty$ and
$$
  \sup_{n\geq 1} \max_{k\leq n} \frac{nb_k^2}{\sigma_n^2} = c < \infty. \tag{2}
$$
Then for for $\frac1{\sigma_n} [S_n - E(S_n| \Cal F_0)]$ a quenched CLT holds true.
\endproclaim

\underbar{Remark 1}: If the sums $b_k$ converge to a limit $b$ such that $\sigma_n^2/n \to b^2$ then the Heyde's condition (cf\. e.g\. \cite{HaHe}) 
is satisfied and we get a CLT.
As proved in \cite{VWo14}, for $S_n - E(S_n| \Cal F_0)$ this CLT is not quenched. Our theorem shows that it is quenched
in the particular case when $(f\circ T^i)$ is a stationary linear process. Recall that any regular stationary process is a sum of 
a countable set of mutually orthogonal stationary linear processes (see \cite{VWoZ}).

\underbar{Remark 2}: For a stationary linear process Theorem 1, in fact, proves more: we get a quenched CLT for $S_n - E(S_n| \Cal F_0)$ 
as soon as $\sum_{k=1}^\infty a_k^2 <\infty$ $\liminf_{n\to\infty} \sigma_n^2/n >0$, and the sequence of $ b_k = \sum_{i=0}^{k} a_i$ is bounded.

\demo{Proof} 
We have to prove a quenched CLT for the triangular array of random variables $b_{n-k}e\circ T^k/\sigma_n$, $k=1\dots,n$, $n=1,2,\dots$.
To do so we use Lachout's refinement \cite{L} of the McLeish's central limit theorem (\cite{Mc}),
applied to regular conditional probabilities with
respect to the $\sigma$-algebra $\Cal F_0$. We thus will prove
\roster
\item"(a)" $E\big(\max_{k\leq n} \frac{|b_{n-k}|}{\sigma_n}|e\circ T^k| \,\big|\, \Cal F_0\big)\to 0$ a.s\.,
\item"(b)" $\sum_{k=1}^n \frac{b_{n-k}^2 e^2\circ T^k}{\sigma_n^2}$ converge to a constant a.s..
\endroster

By (2), 
$$
  \frac{b_{n-k}^2}{\sigma_n^2} \leq \frac{c}{n}
$$
for all $n$, $1\leq k\leq n$, hence (a) follows in the same way as in \cite{VWo14}.

We thus have to prove (b).

Remark that $\sum_{k=1}^n \frac{b_{n-k}^2 e^2\circ T^k}{\sigma_n^2} \to 1$ implies $b_n^2/\sigma_n^2 \to 0$. From this we deduce that
$\max_{k\leq n} \frac{|b_{n-k}|}{\sigma_n} = \max_{k\leq n} \frac{|b_{k}|}{\sigma_n} \to 0$.

Denote 
$$
  T_nf = \frac1{\sigma_n^2} \sum_{k=1}^n b_{n-k}^2 f\circ T^k, \quad f\in L^1.
$$
Recall the Banach's principle: \newline
If
\roster
\item"(i)" $T_n$ : $L^1 \to L^1$ are continuous,
\item"(ii)" for every $f\in L^1$, $\sup_n |T_nf| <\infty$ a.e.,
\item"(iii)" there is a dense subset of $h\in L^1$ for which $T_nh$ converges a.s.,
\endroster
then
for all $f\in L^1$n $T_nf$ converge a.s..

(i) follows from the definition.

For (ii),
$$
   |T_nf| \leq \frac1{\sigma_n^2} \sum_{k=1}^n b_{n-k}^2 |f|\circ T^k \leq \frac{c}{n} \sum_{k=1}^n |f|\circ T^k
$$
hence, by the Birkhoff's ergodic theorem, 
$$
  \sup |T_nf| < \infty\,\,\,\,a.s.\quad \forall f\in L^1.
$$
\medskip

Let us prove (iii). Let $f = g-g\circ T$, $g\in L^\infty$. Then
$$\gathered
  T_n f = \frac1{\sigma_n^2} \sum_{k=1}^n b_{n-k}^2 [g\circ T^k - g\circ T^{k+1}] = \\
  = \frac1{\sigma_n^2} \sum_{k=1}^n [b_{n-k}^2 - b_{n-k+1}^2] g\circ T^k + \frac{b_n^2}{\sigma_n^2} g\circ T \leq \\
  \leq \frac1{\sigma_n^2} \sqrt{\sum_{k=1}^n (b_{n-k} + b_{n-k+1})^2} \sqrt{\sum_{k=1}^n a_{n-k+1}^2} \|g\|_\infty \leq
  \frac2{\sigma_n} A \|g\|_\infty
  \endgathered \tag{3}
$$  
where $A^2 = \sum_{k=1}^\infty a_k^2 <\infty$.

The set of functions $c+g-g\circ T$, $c\in \Bbb R$, $g\in L^\infty$, is dense in $L^1$;
for $f' = g-g\circ T$ we have $T_n f \to 0$ a.s. by the calculation above, for $f''=c$ we have $T_nf''=f''$ hence the convergence
takes place for $f=f'+f''$.
\medskip

We thus conclude that 
$$
  \frac1{\sigma_n^2} \sum_{k=1}^n b_{n-k}^2 e^2\circ T^k \tag{4}
$$
converges almost surely. It remains to prove that the limit is a constant.

Let $f^*$ be the limit in (4). 
Using the same  calculation as in (3) when replacing $g-g\circ T$ by $f-f\circ T$ and $L^\infty$ norm by $L^1$ norm
we prove $f^* = f^*\circ T$. By ergodicity $f^*$ is thus constant.
\enddemo
\qed

\bigskip

\subheading{2. A Non-quenched CLT}

\proclaim{Theorem 2} There exists a stationary linear process $(f\circ T^i)$ with martingale difference innovations such that
\roster
\item"({\bf i})" the Maxwell-Woodroofe condition 
$$
  \sum_{n=1}^\infty \frac{\|E(S_n(f)\,|\,\Cal F_0)\|_2}{n^{3/2}} < \infty
$$  
is satisfied,
\item"({\bf ii})" for $\sigma_n = \|S_n(f)\|_2$, $\sigma_n \to\infty$, $\sigma_n /\sqrt n \to 0$, $\|E(S_n(f)\,|\,\Cal F_0)\|_2/\sigma_n \to 0$, i.e\.
$\|S_n(f) - E(S_n(f)\,|\,\Cal F_0)\|_2/\sigma_n \to 0$,
\item"({\bf iii})" $S_n(f)/\sigma_n$ converge in distribution to the standard normal law $N(0, 1)$,
\item"({\bf iv})" the convergence is not quenched neither for $S_n(f)/\sigma_n$ nor for $(S_n(f) - E(S_n(f)\,|\,\Cal F_0))/\sigma_n$.
\endroster
\endproclaim

\demo{Proof}
We will find a filtration $(\Cal F_i)_i$ such that $\Cal F_{i+1} = T^{-1}\Cal F_i$ and $e\in L^2(\Cal F_0)\ominus L^2(\Cal F_{-1})$,
$\|e\|_2 =1$. The construction of $e$ and $(\Cal F_i)_i$ will be presented later; it is needed for the proof of ({\bf iv}) only.

We define a function $f$ by
$$
  f = e + \sum_{k=1}^{\infty} \frac{-\gamma_k}{V_k} \sum_{i=1}^{V_k} U^{-i}e
$$
where $\gamma_k>0$, $\sum_{k=1}^{\infty} \gamma_k = 1$, $V_k \nearrow \infty$, are such that

$$\gather
  \|S_n(f)\|_2\to \infty, \quad \frac{\|S_n(f)\|_2}{\sqrt n} \to 0, \quad
  \frac{\|E(S_n(f)\,|\,\Cal F_0)\|_2}{\|S_n(f)\|_2} \to 0, \\
  \quad\text{and}\quad
  \frac{S_n(f)}{\|S_n(f)\|_2} \to N(0, 1).
  \endgather
$$

To do so, we define
$$
  \gamma_k = \frac2{k+2} \prod_{j=1}^{k} \big(1 - \frac1{j+1}\big), \quad k=1,2,\dots.
$$
From
$$
  \gamma_k = 2 \prod_{j=1}^{k} \big(1 - \frac1{j+1}\big) - 2 \prod_{j=1}^{k+1} \big(1 - \frac1{j+1}\big)
$$
and
$$
  \lim_{k\to\infty} \prod_{j=1}^{k} \big(1 - \frac1{j+1}\big) = 0
$$
we deduce
$$
   \sum_{j=k}^\infty \gamma_j = 2 \prod_{j=1}^{k} \big(1 - \frac1{j+1}\big).
$$
Therefore,
$$
  \sum_{k=1}^\infty \gamma_k =1,\quad 1 - \sum_{j=1}^{k-1} \gamma_j = 2 \prod_{j=1}^{k} \big(1 - \frac1{j+1}\big) = (k+2)\gamma_k.
$$

The numbers $V_k$ will be specified later.

We have
$$
  \big\|\frac{\gamma_k}{V_k} \sum_{i=1}^{V_k} U^{-i}e \big\|_2 = \frac{\gamma_k}{\sqrt{V_k}}
$$
which guarantees that
$$
  f = \sum_{k=1}^\infty \gamma_k \Big( e - \frac{1}{V_k} \sum_{i=1}^{V_k} U^{-i}e \Big) =
  \sum_{k=1}^\infty \gamma_k f_k \in L^2.
$$
For every $k$, $f_k = e - \frac{1}{V_k} \sum_{i=1}^{V_k} U^{-i}e$ is a coboundary $f_k = g_k - Ug_k$ where
$$
  g_k = -\frac1{V_k} \sum_{j=1}^{V_k} jU^{-V_k-1+j} e.
$$  
For $h= \sum_{i=0}^\infty c_iU^{-i}e$ we have
$$
  S_n(h) = \sum_{j=0}^{n-1} \sum_{i=0}^\infty c_i U^{j-i}e = \sum_{u=-\infty}^{n-1} \sum_{j=\max\{0, u\}}^{n-1} c_{j-u}U^ue. \tag{5}
$$
Let us put
$$
  f_k =  e - \frac{1}{V_k} \sum_{i=1}^{V_k} U^{-i}e = \sum_{i=0}^\infty c_{k,i}U^{-i}e.
$$
We then have
$$\gathered
  \big| \sum_{j=\max\{0, u\}}^{n-1} c_{k,j-u}\big| \leq 1\,\,\,\,\text{for every}\,\,\,\, u, \\
  \big| \sum_{j=\max\{0, u\}}^{n-1} c_{k,j-u}\big| \leq \frac{n}{V_k} \,\,\,\,\text{for}\,\,\,\,-1\geq u \geq -V_k, \\
  \big| \sum_{j=\max\{0, u\}}^{n-1} c_{k,j-u}\big| = 0 \,\,\,\,\text{for}\,\,\,\,u<-V_k.
  \endgathered \tag6
$$  
We deduce that for every $V_k$,
$$
  \big\|S_n\big( e - \frac{1}{V_k} \sum_{i=1}^{V_k} U^{-i}e \big)\big\|_2 \leq \sqrt{2 n}
$$
hence
$$
  \frac{\|S_n(f)\|_2}{\sqrt n} \leq \sum_{k=1}^\infty \gamma_k \frac{\|S_n(f_k)\|_2}{\sqrt n} \to 0 \tag7
$$
by the Lebesgue Dominated Convergence Theorem.

\noindent For $V_k \leq n <V_{k+1}$ we using (5), (6) get
$$
  \|S_n(f)\|_2 \geq \|S_n(f) - E(S_n(f)\,|\,\Cal F_0)\|_2 \geq \sqrt n \big(1 - \sum_{j=1}^{k+1} \gamma_j -
  \sum_{j=k+2}^{\infty} \frac{V_{k+1}}{V_j}\gamma_j \big). \tag8
$$
Recall that
$$
  1 - \sum_{j=1}^{k+1} \gamma_j = \sum_{j=k+2}^\infty \gamma_j
$$
hence for $V_j$ growing exponentially fast and
$$
  \sqrt{V_k} \big(1 - \sum_{j=1}^{k+1} \gamma_j \big) \to \infty
$$
we get
$$
  \|S_n(f)||_2 \to \infty. \tag9
$$
For $V_k$ growing sufficiently fast we have, for $V_k \leq n <V_{k+1}$,
$$
  \|E(S_n(f)\,|\,\Cal F_0)\|_2 \sim (\gamma_k+\gamma_{k+1})\sqrt n. \tag{10}
$$
Indeed, $f_j = g_-Ug_j$ where $g_j\in L^2$ hence $\|E(S_n(f)\,|\,\Cal F_0)\|_2 =o(n)$ for $j\leq k-1$ and if $V_j$ are large enough
for $j\geq k+2$, $\|E(S_n(f)\,|\,\Cal F_0)\|_2 =o(n)$ as well.\newline
By definition,
$$
  \gamma_k = \frac1{k+2} \big(1 - \sum_{j=1}^{k-1} \gamma_j \big)
$$
hence by (8) and (10) we have
$$
  \frac{\|E(S_n(f)\,|\,\Cal F_0)\|_2}{\|S_n(f)\|_2} \to 0. \tag{11}
$$
From (9), (7), and (11) we get ({\bf ii}).
\medskip

By \cite{PU06} and $\sigma_n^2 \to \infty$ we get the central limit theorem ({\bf iii}).
\bigskip

\centerline{\it The Maxwell-Woodroofe condition}
\medskip

Denote $h_k = \frac{-\gamma_k}{V_k} \sum_{i=1}^{V_k} U^{-i}e$ and define $f' = \sum_{k=1}^\infty h_k$; we thus have $f=e-f'$.
By (6), for any $k\geq 1$, 
$$
  \|S_n(h_k)\|_2^2 \leq \gamma_k^2 2V_k \Big(\frac{n}{V_k}\Big)^2 = 2 \gamma_k^2 \frac{n^2}{V_k}, \quad n=1,2,...,V_k
$$
and
$$
  \|E(S_n(h_k)\,|\,\Cal F_0)\|_2 = \|E(S_{V_k}(h_k)\,|\,\Cal F_0)\|_2 \leq \gamma_k \sqrt{V_k} , \quad n\geq V_k
$$
hence
$$
  \sum_{n=1}^\infty \frac{\|E(S_n(h_k)\,|\,\Cal F_0)\|_2}{n^{3/2}} \leq \sqrt 2 \gamma_k \frac{1}{\sqrt{V_k}} 
  \sum_{n=1}^{V_k} \frac{1}{\sqrt{n}} +
  \gamma_k \sqrt{V_k} \sum_{n=V_k+1}^\infty \frac1{n^{3/2}} \leq C\gamma_k
$$
for some constant $C$. Therefore,
$$
  \sum_{n=1}^\infty \frac{\|E(S_n(f')\,|\,\Cal F_0)\|_2}{n^{3/2}} \leq
  \sum_{k=1}^\infty \sum_{n=1}^\infty \frac{\|E(S_n(h_k)\,|\,\Cal F_0)\|_2}{n^{3/2}}
  \leq C \sum_{k=1}^\infty \gamma_k = C <\infty
$$
and ({\bf i}) follows.

\bigskip

\centerline{\it The filtration and $e$}
\medskip

Let $\Cal B_k',\Cal B_l'' \subset \Cal A$, $k, l=1,2,\dots$, be mutually independent $\sigma$-algebras,
$$\gather
  \Cal B_k'\subset T^{-1}\Cal B_k',\quad \Cal B_l''\subset T^{-1}\Cal B_l'',\\
  \cap_{j=1}^\infty T^{j}\Cal B_k' = \{\Omega, \emptyset\},\quad  \cap_{j=1}^\infty T^{j}\Cal B_l'' = \{\Omega, \emptyset\}
  \endgather
$$
(modulo sets of measure 0 or 1) for every $k, l$. 

$\xi_k\circ T^i$ are iid $\Cal B_k'$-measurable random variables, $\mu(\xi_k=1) = 1/2 = \mu(\xi_k=-1)$ for all $i$.

All of this can be constructed by taking finite alphabets $\Bbb A_k'$ and $\Bbb A_l''$, $k,l=1,2,\dots$,
$\Omega_k' = \underset i\in \Bbb Z \to{\times} \Bbb A_{k,i}'$ where $\Bbb A_{k,i}'$ are identical copies of $\Bbb A_k'$, similarly we define
$\Omega_l''$, $k,l=1,2,\dots$. On the sets $\Omega_k'$ and $\Omega_l''$ we define product $\sigma$-algebras, product measures, and 
left shift transformations $T_k'$, $T_l''$. $\Omega$ is the product of all $\Omega_k'$ and $\Omega_l''$ with the product $\sigma$-algebra 
$\Cal A$,
product (probability) measure $\mu$, and product transformation $T$. For projections $\xi_k$ and $\zeta_l$ of $\Omega$ onto $\Bbb A_{k,0}'$ 
and $\Bbb A_{l,0}''$ we thus get mutually independent processes of iid $(\xi_k\circ T^i)_i$, $(\zeta_l\circ T^i)_i$. We suppose that
$A_k = \{-1, 1\}$, $k=1,2, \dots$, and $\mu(\xi_k=1) = 1/2 = \mu(\xi_k=-1)$. For $\Cal B_k'$ we take the past $\sigma$-algebras
$\sigma\{\xi_k\circ T^i : i\leq 0\}$ and for $\Cal B_l''$ we take the past $\sigma$-algebras $\sigma\{\zeta_l\circ T^i : i\leq 0\}$.
The properties above can be easily verified, the latter follow from Kolmogorov's 0-1 law.

We thus have 
that $\xi_k\circ T^i$ are iid $\Cal B_k'$-measurable random variables, $\mu(\xi_k=1) = 1/2 = \mu(\xi_k=-1)$.
For $k=1,2,\dots$, let $A_k\in \Cal B_k''$ be sets and $N_k$ positive integers such that $T^{-i}A_k$, $i=0,\dots,3N_k$ are mutually disjoint 
(hence $\{T^{-i}A_k :\, i=0,\dots,3N_k\}$ are Rokhlin towers) and $\mu(A_k)= \frac1{4N_k}$;
the values of $N_k$ will be specified later (existence of Rokhlin Towers is proved e.g\. in \cite{CSF}). We suppose that
$$
  \sum_{k=1}^\infty \mu(A_k) < \frac12. \tag{12}
$$

By $\Cal B''$ we define the $\sigma$-algebra generated by all $T^{-i}\Cal B_k''$, $i\in\Bbb Z$, $k=1,2,\dots$; we thus have $T^{-1}\Cal B''
=\Cal B''$, all Rokhlin towers defined above are $\Cal B''$-measurable.\newline
By $\Cal F_j$ we denote the $\sigma$-algebra generated by $\Cal B''$ and all $\xi_k\circ T^i$, $i\leq j$, $k=1,2,\dots$; notice that
$T^{-1}\Cal F_j = \Cal F_{j+1}$.

For $c = 2 (\sum_{k=1}^\infty \frac1{k^3})^{-1/2}$ we define
$$
  e_k = c\xi_k\frac{\sqrt{N_k}}{k^{3/2}}1_{A_k}, \quad e = \sum_{k=1}^\infty e_k.
$$
Notice that $\|e_k\|_2 = \frac{c}{2 k^{3/2}}$ hence $e\in L^2$. By definition, $e$ is $\Cal F_0$-measurable. 
By definition, $e_k$ are mutually independent hence $\|e\|_2^2 = \frac{c^2}4 \sum_{k=1}^\infty \frac1{k^3}$; we thus have 
$\|e\|_2 = 1$. \newline
Because $A_k\in \Cal B''$ and $\xi_k$ is independent of $\Cal F_{-1}$, we have $E(e_k\,|\,\Cal F_{-1}) =0$ for every $k$ hence 
$E(e\,|\,\Cal F_{-1}) =0$, $(U^ie)_i$ is thus a martingale difference sequence adapted to the filtration $(\Cal F_i)$.

Recall that
$$
  f = e + \sum_{k=1}^{\infty} \frac{-\gamma_k}{V_k} \sum_{i=1}^{V_k} U^{-i}e = a_0e - \sum_{i=1}^\infty a_i U^{-i}e
$$
where $a_0=1$, $a_i>0$ for all $i\geq 1$, and $\sum_{i=1}^\infty a_i = 1$.

By $m_{\omega}$ we will denote regular conditional probabilities w.r.t\. $\Cal F_0$ ($\Cal A$ is a Borel $\sigma$-algebra of a Polish space
hence the regular conditional probabilities exist). Notice that all sets $T^{-i}A_k$, $k=1,2,\dots$, 
$i\in\Bbb Z$, are $\Cal F_0$-measurable hence $m_{\omega}(T^{-i}A_k) = 0$ (if $\omega\not\in T^{-i}A_k$) or $m_{\omega}(T^{-i}A_k) = 1$ 
(if $\omega\in T^{-i}A_k$). 
\smallskip

Let us fix a $k\geq 1$ such that
$$
  \sum_{i=N_k}^\infty a_i < \frac12 
$$
and denote 
$$
  A_k' = A_k \setminus \underset j {\neq k} \to{\bigcup} A_j.
$$
By (12) and independence, $\mu(A_k') \geq \frac12 \mu(A_k)$.

Let us denote 
$$$\Bbb Z^d$ action
  A_k' = \big(A_k \setminus \underset j {\neq k} \to{\bigcup} A_j\big).
$$  
The sets $A_k',\dots, T^{-3N_k+1}A_k'$ are mutually disjoint and $\mu(A_k)\geq \frac1{4N_k}$ hence
$$
  \mu\Big(\underset N=0 \to{\overset N_k-1\to{\bigcup}} T^{-N+1} A_k'\Big) \geq \frac1{8}.\tag{13}
%  \big(A_k \setminus \underset j\neq k \to{\cup} A_j\big) \Big) \geq \frac1{16}.\tag{13}
$$  
Let $N_k\leq N < N_{k+1}$.
For $\omega \in T^{-N+1} A_k'$
%\Big(A_k \setminus \underset j \neq k \to{\bigcup} A_j\Big)$
we have
$$
  |U^{N-1}e(\omega)| = |U^{N-1}e_k(\omega)| = c\frac{\sqrt{N_k}}{k^{3/2}} 1_{T^{-N+1}A_k}(\omega) = c\frac{\sqrt{N_k}}{k^{3/2}}
$$
and
$$\multline
%  S_N(f) = \sum_{j=0}^{N-1} \sum_{i=0}^\infty a'_i U^{j-i}e = \sum_{u=-\infty}^{N-1} \sum_{j=\max\{0, u\}}^{N-1} a'_{j-u}U^ue
  S_N(f) = \sum_{j=0}^{N-1} \sum_{i=0}^\infty a_i U^{j-i}e = \sum_{u=-\infty}^{N-1} \sum_{j=\max\{0, u\}}^{N-1} a_{j-u}U^ue =   \\
  = U^{N-1}e_k
  - c \frac{\sqrt{N_k}}{k^{3/2}} \sum_{u=-\infty}^{0} \sum_{j=\max\{0, u\}}^{N-1} a_{j-u} 1_{\{T^{-u}A_k\}} U^u\xi_k - \\
  - c \frac{\sqrt{N_k}}{k^{3/2}} \sum_{u=1}^{N-2} \sum_{j=\max\{0, u\}}^{N-1} a_{j-u} 1_{\{T^{-u}A_k\}} $\Bbb Z^d$ actionU^u\xi_k - \\
  - c \sum_{l\neq k} \frac{\sqrt{N_l}}{l^{3/2}} \sum_{u=-\infty}^{0} \sum_{j=\max\{0, u\}}^{N-1} a_{j-u}  1_{\{T^{-u}A_l\}}U^u\xi_l -\\
  - c \sum_{l\neq k} \frac{\sqrt{N_l}}{l^{3/2}} \sum_{u=1}^{N-2} \sum_{j=\max\{0, u\}}^{N-1} a_{j-u}  1_{\{T^{-u}A_l\}}U^u\xi_l = \\
  = U^{N-1}e_k - I - II - III-IV.
  \endmultline \tag{14}
$$
We have 
$$
  m_\omega\big(U^{N-1}e(\omega) = c\frac{\sqrt{N_k}}{k^{3/2}}\big) = m_\omega\big(U^{N-1}e(\omega) = -c\frac{\sqrt{N_k}}{k^{3/2}}\big) 
  = \frac12
$$
and $U^{N-1}e_k$ is independent of $I$, $II$, $III$, $IV$ (with respect to $m_\omega$) hence
$$
  m_\omega\big( |S_N(f)| \geq c\frac{\sqrt{N_k}}{k^{3/2}}\big) \geq \frac14. \tag{15}
$$  
Recall that 
$$
  \sigma_n/\sqrt n = \|S_n(f)\|_2 /\sqrt n \to 0.
$$
We thus can choose $N_k$, $N_{k+1}$ such that $2^k\sigma_{N_k} \leq c \frac{\sqrt{N_k}}{k^{3/2}}$
and $N_{k+1} = 2N_k$ for $k$ odd. 
From this and (15) it follows that there exists a $d>0$ such that for $k$ odd and sufficiently big, the Prokhorov (metric) distance  between 
the distribution of $S_N(f)/\sigma_N$ and $\Cal N(0, 1)$ is bigger than $d$.\newline
For the given odd $k$ there thus exists a set $B_k$ of measure bigger than $1/8$ (cf\. (13)) such that for $\omega\in B$ there exists
an $N_k\leq N\leq N_{k+1}$ for which the Prokhorov (metric) distance  between the distribution of $S_N(f)/\sigma_N$ and $\Cal N(0, 1)$ 
is bigger than $d$.

There exists a set $B$ of positive measure such that for $\omega\in B$ there is an infinite sequence of
$k$ (odd) and $N_k\leq N\leq N_{k+1}$ such that for the probability $m_\omega$, the laws of $S_N(f)/\sigma_N$ do not weakly converge to
$\Cal N(0, 1)$ hence the CLT for $S_n/\sigma_n$ is not quenched.
\smallskip

In (14) we can note that $E(S_N | \Cal F_0) = -I-III$ and we deduce in the same way as above that the CLT for $S_n - E(S_n | \Cal F_0)$
is not quenched either. This proves ({\bf iv}).

\comment
\bigskip

By the definition of the Rokhlin towers, from $1_{\{T^{-u}A_k\}}(\omega) =1$ it follows 
$1_{\{T^{-v} A_k\}}(\omega) =0$ for $|u-v| \leq 3N_k$; we thus have
$$\gather
  1_{T^{-N+1}A_k}\big|\sum_{u=-\infty}^{N-2} \sum_{j=\max\{0, u\}}^{N-1} a_{j-u} 1_{\{T^{-u}A_k\}}\big| = \\
  = 1_{T^{-N+1}A_k}\big|\sum_{u=-\infty}^{N-3N_k-1} \sum_{j=\max\{0, u\}}^{N-1} a_{j-u} 1_{\{T^{-u}A_k\}}\big| \leq 
  1_{T^{-N+1}A_k} \sum_{i=N_k}^\infty a_i,
  \endgather
$$
by (14) thus 
$$
  | I | \leq  1_{T^{-N+1}A_k}\big|\sum_{u=-\infty}^{N-2} \sum_{j=\max\{0, u\}}^{N-1} a_{j-u} 1_{\{T^{-u}A_k\}}\big| < \frac12.
$$  
The sum $II$ is $\Cal F_0$-measurable.
The random variables $U^u\xi_l$, $u\geq 1$, are conditionally independent of $U^{N-1}e_k$ with respect to $\Cal F_0$.

Let $\omega \in T^{-N+1} A_k'$ and $m_\omega$ be the regular conditional probability with respect to $\Cal F_0$.
$S_N(f)$ is a sum of independent (for the measure $m_\omega$) random variables $U^{N-1}e_k$, $I$, and $III$ and an $m_\omega$-a.s\.
constant $II$. By the previous calculation we have
$$
  |U^{N-1}e_k - I| \geq \frac{c}2 \frac{\sqrt{N_k}}{k^{3/2}}
$$
($m_\omega$-a.s\.). \newline
Recall that 
$$
  \sigma_n/\sqrt n = \|S_n(f)\|_2 /\sqrt n \to 0.
$$
We thus can choose $N_k$, $N_{k+1}$ such that $\sigma_N \leq \frac{c}4 \frac{\sqrt{N_k}}{k^{3/2}}$

we have 
$$
  E_{m_\omega}(S_N(f) - E_{m_\omega}S_N(f))^2 \geq \frac{c^2}4 \frac{N_k}{k^3}.
$$

\medskip

Recall that 
$$
  \sigma_n/\sqrt n = \|S_n(f)\|_2 /\sqrt n \to 0.
$$
We can choose $N_k$ such that for $N_k\leq n < N_{k+1}$ it is $\frac{\sqrt{N_k}}{k^{3/2}} \geq 2^k \sigma_n$ ($\sigma_n = \|S_n(f)\|_2$).
We then get
$$
  |E_{m_\omega}S_N(f)| \geq c2^{k-1} \sigma_n.
$$
\medskip

The sets $B_k$ are each of probability greater than $1/32$. There thus exists a set $B\in \Cal F_0$, $\mu(B)>0$, such that
for $\omega\in B$ there exist infinitely many $k$ and $N_{k-1}\leq N\leq N_k$ such that
$$
  |E_{m_\omega}S_N(f)| \geq c2^{k-1} \sigma_n.
$$
and the central limit theorem ($\gamma$) is thus not quenched.

\medskip
\endcomment

\enddemo
\qed

\underbar{Remark}. We have proved that there exists a $C>0$ such that 
$$
  \mu(\max_{1\leq n\leq N_{k+1}} |S_n(f)| > C2^{k} \sigma_{N_{k+1}}) > \frac1{32}
$$
for infinitely many $k$ hence the Donsker invariance principle is not satisfied (more precisely, by \cite{PU05} it is satisfied 
for $S_n(f)/\sqrt n$ with a convergence towards a degenerate limit, but it is not satisfied for $S_n(f)/\sigma_n$).

\enddocument
\end

\ref \key B \by Billingsley, P. \paper The Lindeberg-L\'evy theorem for martingales \jour Proc. Amer. Math. Soc.
\vol 12 \pages 788-792 \yr 1961 \endref 
\ref \key C-F-S \by Cornfeld, I.P., Fomin, S.V., and Sinai, Ya.G. \book Ergodic Theory \publ Springer \publaddr Berlin \yr 1982 \endref
\ref \key Cu 1 \by Cuny, Ch. \paper Pointwise ergodic theorems with rate and application to limit theorems for stationary processes
\paperinfo sumitted for publication, arXiv:0904.0185 \yr 2009 \endref
\ref \key Cu 2 \by Cuny, Ch. \paper Norm convergence of some power-series of operators in $L^p$ with applications in ergodic theory
\paperinfo sumitted for publication \yr 2009 \endref
\ref \key D-L 1 \by Derriennic, Y. and Lin, M. \paper Sur le th\'eor\`eme limite central de Kipnis et Varadhan pour les cha\^\i nes 
r\'eversibles ou normales \jour CRAS 323 \pages 1053-1057 \yr 1996 \endref 
\ref \key D-L 2 \by Derriennic, Y and Lin, M. \paper The central limit theorem for Markov
chains with normal transition operators, started at a point \jour Probab. Theory Relat. Fields
\vol 119 \pages 509-528 \yr 2001 \endref
\ref \key G \by Gordin, M.I. \paper A central limit theorem for stationary processes \jour Soviet Math. Dokl. \vol 10 
\pages 1174-1176 \yr 1969 \endref
\ref \key G-Ho \by Gordin, M.I. and Holzmann, H. \paper The central limit theorem for stationary Markov chains under invariant splittings 
\jour Stochastics and Dynamics \vol 4 \pages 15-30 \yr 2004 \endref 
\ref \key G-L 1 \by Gordin, M.I. and Lif\v sic, B.A. \paper Central limit theorem for statioanry processes \jour Soviet Math. Doklady 19 
\pages 392-394 \yr 1978 \endref 
\ref \key G-L 2 \by Gordin, M.I. and Lif\v sic, B.A. \paper A remark about a Markov process with normal transition operator
\paperinfo In: Third Vilnius Conference on Probability and Statistics \vol 1 \pages 147-148 \yr 1981 \endref
\ref \key G-L 3 \by Gordin, M.I. and Lif\v sic, B.A. \paper The central limit theorem for Markov processes with normal transition operator,
and a strong form of the central limit theorem \paperinfo §IV.7 and §IV.8 in Borodin and Ibragimov, Limit theorems for functionals of
random walks, Proc. Steklov Inst. Math. 195(1994), English translation AMS (1995) \endref
\ref \key Ha-He \by Hall, P. and Heyde, C.C. \book Martingale Limit Theory           
and its Application \publ Academic Press \publaddr New York \yr 1980 \endref 
\ref \key I \by Ibragimov, I.A.  \paper A central limit theorem for a class of dependent random variables 
\jour Theory Probab. Appl. \vol 8 \pages 83-89 \yr 1963 \endref
\ref \key K-V \by Kipnis, C. and Varadhan, S.R.S. \paper Central limit theorem for additive functionals of reversible
Markov processes and applications to simple exclusions \jour Comm. Math. Phys. \vol 104 \pages 1-19 \yr 1986 \endref
\ref \key Kl-Vo 1 \by Klicnarov\'a, J. and Voln\'y, D. \paper An invariance principle for non adapted processes
\jour C.R. Acad. Sci. Paris Ser 1 \vol 345/5 \pages 283-287 \yr 2007 \endref
\ref \key Kl-Vo 2 \by Klicnarov\'a, J. and Voln\'y, D. \paper Exactness of a Wu-Woodroofe's approximation with linear growth of variances \jour Stoch. Proc. and their Appl. \vol 119 \pages 2158-2165 \yr 2009 \endref
\ref \key M-Wo \by Maxwell, M. and Woodroofe, M. \paper Central limit theorems for additive
functionals of Markov chains \jour Ann. Probab. \vol 28 \pages 713-724 \yr 2000 \endref
\ref \key P-U \by Peligrad, M. and Utev, S. \paper A new maximal inequality
and invariance principle for stationary sequences \jour Ann. Probab.
\vol 33 \pages 798-815 \yr 2005 \endref
\ref \key R \by Rosenblatt, M. \book Markov Processes: Structure and asymptotic behavior \publ Springer \publaddr Berlin \yr 1971 
\endref
\ref \key Vo 1 \by Voln\'y, D. \paper Approximating martingales and the
central limit theorem for strictly stationary processes
\jour Stochastic Processes and their Applications \vol 44 \pages 41-74
\yr 1993 \endref
\ref \key Vo 2 \by Voln\'y, D. \paper Martingale approximation of non adapted stochastic processes
with nonlinear growth of variance \paperinfo Dependence in Probability and Statistics
Series: Lecture Notes in Statistics, Vol. 187 Bertail, Patrice; Doukhan, Paul; Soulier, Philippe (Eds.) 
\yr 2006 \endref
\ref \key Wo \by Woodroofe, M. \paper A central limit theorem for functions of a Markov chain with applications to
shifts \jour Stoch. Proc. and their Appl. \vol 41 \pages 31-42 \yr 1992 \endref
\ref  \key Wu-Wo \by Wu, W.B. and Woodroofe, M. \paper Martingale approximation for 
sums of stationary processes \jour Ann. Probab. \vol 32 \pages 1674-1690 \yr 2004 \endref  
\endRefs
\enddocument
\end

\end

drrrrrrrrrrrrrrrrrrrrrrrrrrrrkgfvvvvvvvvvvvvvve central limit theorem for Markov processes with normal transition operator,
and a strong form of the central limit theorem \paperinfo §IV.7 and §IV.8 in Borodin and Ibragimov, Limit theorems for functionals of
random walks, Proc. Steklov Inst. Math. 195(1994), English translation AMS (1995) \endref
\ref \key Ha-He \by Hall, P. and Heyde, C.C. \book Martingale Limit Theory           
and its Application \publ Academic Press \publaddr New York \yr 1980 \endref 
\ref \key I \by Ibragimov, I.A.  \paper A central limit theorem for a class of dependent random variables 
\jour Theory Probab. Appl. \vol 8 \pages 83-89 \yr 1963 \endref
\ref \key K-V \by Kipnis, C. and Varadhan, S.R.S. \paper Central limit theorem for additive functionals of reversible
Markov processes and applications to simple exclusions \jour Comm. Math. Phys. \vol 104 \pages 1-19 \yr 1986 \endref
\ref \key Kl-Vo 1 \by Klicnarov\'a, J. and Voln\'y, D. \paper An invariance principle for non adapted processes
\jour C.R. Acad. Sci. Paris Ser 1 \vol 345/5 \pages 283-287 \yr 2007 \endref
\ref \key Kl-Vo 2 \by Klicnarov\'a, J. and Voln\'y, D. \paper Exactness of a Wu-Woodroofe's approximation with linear growth of variances \jour Stoch. Proc. and their Appl. \vol 119 \pages 2158-2165 \yr 2009 \endref
\ref \key M-Wo \by Maxwell, M. and Woodroofe, M. \paper Central limit theorems for additive
functionals of Markov chains \jour Ann. Probab. \vol 28 \pages 713-724 \yr 2000 \endref
\ref \key P-U \by Peligrad, M. and Utev, S. \paper A new maximal inequality
and invariance principle for stationary sequences \jour Ann. Probab.
\vol 33 \pages 798-815 \yr 2005 \endref
\ref \key R \by Rosenblatt, M. \book Markov Processes: Structure and asymptotic behavior \publ Springer \publaddr Berlin \yr 1971 
\endref
\ref \key Vo 1 \by Voln\'y, D. \paper Approximating martingales and the
central limit theorem for strictly stationary processes
\jour Stochastic Processes and their Applications \vol 44 \pages 41-74
\yr 1993 \endref
\ref \key Vo 2 \by Voln\'y, D. \paper Martingale approximation of non adapted stochastic processes
with nonlinear growth of variance \paperinfo Dependence in Probability and Statistics
Series: Lecture Notes in Statistics, Vol. 187 Bertail, Patrice; Doukhan, Paul; Soulier, Philippe (Eds.) 
\yr 2006 \endref
\ref \key Wo \by Woodroofe, M. \paper A central limit theorem for functions of a Markov chain with applications to
shifts \jour Stoch. Proc. and their Appl. \vol 41 \pages 31-42 \yr 1992 \endref
\ref  \key Wu-Wo \by Wu, W.B. and Woodroofe, M. \paper Martingale approximation for 
sums of stationary processes \jour Ann. Probab. \vol 32 \pages 1674-1690 \yr 2004 \endref  
\endRefs

\Refs
\widestnumber\key{VWo14}

\ref \key CFS \by Cornfeld, I.P., Fomin, S.V., and Sinai, Ya.G. \book Ergodic Theory \publ Springer \publaddr Berlin \yr 1982 \endref

\ref \key CuMe \by C. Cuny, F. Merlev\`ede \paper On martingale approximation and the quenched weak invariance principle
\jour Ann. Probab. \vol 42 \pages 760-793 \yr 2014 \endref

\ref \key CuP \by C. Cuny, M. Peligrad \paper Central limit theorem started at a point for stationary processes and additive functinals 
of reversible Markov chains
\jour J. of Theoretical Probability \vol 25 \pages 171-188 \yr 2012 \endref

\ref \key CuV \by C. Cuny, D. Voln\'y \paper A quenched invariance principle for stationary processes 
\jour ALEA Lat. Am. J. Probab. Math. Stat. \vol 10 \pages 107-115 \yr 2013 \endref

\ref \key DLi \by Y. Derriennic, M. Lin \paper The central limit theorem for Markov chains with normal transition operators, started at a point 
\jour Probab. Theory Relat. Fields \vol 119 \pages 509-528 \yr 2001 \endref

\ref \key HaHe \by Hall, P. and Heyde, C.C. \book Martingale Limit Theory           
and its Application \publ Academic Press \publaddr New York \yr 1980 \endref 

\ref \key K \by U. Krengel \book Ergodic Theorems 
\publ Walter de Gruyter \publaddr Berlin, New York \yr 1985 \endref

\ref \key L \by P. Lachout \paper A note on the martingale central limit theorem \jour Commentationes Mathematicae Universitatis Carolinae
\vol 26 \pages 637-640 \yr 1985 \endref

\ref \key MWo \by Maxwell, M. and Woodroofe, M. \paper Central limit theorems for additive
functionals of Markov chains \jour Ann. Probab. \vol 28 \pages 713-724 \yr 2000 \endref

\ref \key Mc \by D.L. McLeish \paper Dependent central limit theorems and invariance principles
\jour Ann. Probab. \vol 2 \pages 620-628 \yr 1974 \endref

\ref \key PU05 \by M. Peligrad, S. Utev \paper A new maximal inequality and invariance principle for stationary sequences
\jour Ann. Probab. \vol 33 \pages 798-815 \yr 2005 \endref

\ref \key PU06 \by M. Peligrad, S. Utev \paper Central limit theorem for stationary linear processes
\jour Ann. Probab. \vol 34 \pages 1241-1643 \yr 2006 \endref

\ref \key V \by D. Voln\'y \paper Martingale approximation and optimality of some conditions for the central limit theorem, 
\jour J. of Theoretical Probability \vol 23 \pages 888-903 \yr 2010 \endref

\ref \key VWo10 \by D. Voln\'y, M. Woodroofe\paper An example of non-quenched convergence in the conditional central limit theorem for partial 
sums of a linear process 
\paperinfo  in Dependence in Probability, Analysis and Number Theory. 
A volume in memory of Walter Philipp. Edited by Istvan Berkes, Richard C. Bradley, Herold Dehling, Magda Peligrad and Robert Tichy, 
Kendrick Press, \pages 317-323 \yr 2010 \endref

\ref \key VWoZ \by D. Voln\'y, M. Woodroofe, Ou Zhao \paper Central limit theorems for superlinear processes
\jour Dynamics and Stochastics \vol 11 \pages 71-80 \yr 2011 \endref

\ref \key VWo14 \by D. Voln\'y, M. Woodroofe:  \paper Quenched central limit theorems for sums of stationary processes,
\jour Stat. and Probab. Letters \vol 85 \pages 161-167 \yr 2014 \endref

\endRefs

\Refs
\widestnumber\key{VWo14}

\ref \key CFS \by Cornfeld, I.P., Fomin, S.V., and Sinai, Ya.G. \book Ergodic Theory \publ Springer \publaddr Berlin \yr 1982 \endref

\ref \key CuMe \by C. Cuny, F. Merlev\`ede \paper On martingale approximation and the quenched weak invariance principle
\paperinfo submitted for publication, arXiv:1202.2964[math.PR]\endref

\ref \key CuP \by C. Cuny, M. Peligrad \paper Central limit theorem started at a point for stationary processes and additive functinals 
of reversible Markov chains
\jour J. of Theoreticla Probability \vol 25 \pages 171-188 \yr 2012 \endref

\ref \key DLi \by Y. Derriennic, M. Lin \paper The central limit theorem for Markov chains with normal transition operators, started at a point 
\jour Probab. Theory Relat. Fields \vol 119 \pages 509-528 \yr 2001 \endref

\ref \key HaHe \by Hall, P. and Heyde, C.C. \book Martingale Limit Theory           
and its Application \publ Academic Press \publaddr New York \yr 1980 \endref 

\ref \key K \by U. Krengel \book Ergodic Theorems 
\publ Walter de Gruyter \publaddr Berlin, New York \yr 1985 \endref

\ref \key L \by P. Lachout \paper A note on the martingale central limit theorem \jour Commentationes Mathematicae Universitatis Carolinae
\vol 26 \pages 637-640 \yr 1985 \endref

\ref \key MWo \by Maxwell, M. and Woodroofe, M. \paper Central limit theorems for additive
functionals of Markov chains \jour Ann. Probab. \vol 28 \pages 713-724 \yr 2000 \endref

\ref \key Mc \by D.L. McLeish \paper Dependent central limit theorems and invariance principles
\jour Ann. Probab. \vol 2 \pages 620-628 \yr 1974 \endref

\ref \key PU05 \by M. Peligrad, S. Utev \paper A new maximal inequality and invariance principle for stationary sequences
\jour Ann. Probab. \vol 33 \pages 798-815 \yr 2005 \endref

\ref \key PU06 \by M. Peligrad, S. Utev \paper Central limit theorem for stationary linear processes
\jour Ann. Probab. \vol 34 \pages 1241-1643 \yr 2006 \endref

\ref \key V \by D. Voln\'y \paper Martingale approximation and optimality of some conditions for the central limit theorem, 
\jour J. of Theoretical Probability \vol 23 \pages 888-903 \yr 2010 \endref

\ref \key VWo10 \by D. Voln\'y, M. Woodroofe\paper An example of non-quenched convergence in the conditional central limit theorem for partial 
sums of a linear process 
\paperinfo  in Dependence in Probability, Analysis and Number Theory. 
A volume in memory of Walter Philipp. Edited by Istvan Berkes, Richard C. Bradley, Herold Dehling, Magda Peligrad and Robert Tichy, 
Kendrick Press, \pages 317-323 \yr 2010 \endref

\ref \key VWoZ \by D. Voln\'y, M. Woodroofe, Ou Zhao \paper Central limit theorems for superlinear processes
\jour Dynamics and Stochastics \vol 11 \pages 71-80 \yr 2011 \endref

\ref \key VWo14 \by D. Voln\'y, M. Woodroofe:  \paper Quenched central limit theorems for sums of stationary processes,
\jour Stat. and Probab. Letters \vol 85 \pages 161-167 \yr 2014 \endref

\endRefs
\enddocument
\end

\end